\documentclass[12pt]{amsart}
\usepackage{amsfonts, amssymb}
\usepackage{parskip, fullpage, verbatim}
\usepackage[colorlinks=true]{hyperref}
\usepackage{breakurl}

\newcommand\CA{{\mathcal A}}

\newcommand\CIF{{\mathcal {IF}}} 
 
\newcommand\CIFAC{{\mathcal {IF\!AC}}}

\newcommand\R{{\varrho}}

\newcommand\BBC{{\mathbb C}}
\newcommand\BBK{{\mathbb K}}

\newcommand\codim{\operatorname{codim}}

\newcommand\Der{{\operatorname{Der}}}

\newcommand\pdeg{\operatorname{pdeg}}

\numberwithin{equation}{section}

\theoremstyle{plain}

\newtheorem{theorem}[equation]{Theorem}

\newtheorem{proposition}[equation]{Proposition}
\theoremstyle{definition}
\newtheorem{defn}[equation]{Definition}
\newtheorem{remark}[equation]{Remark}
\newtheorem{example}[equation]{Example}


\begin{document}

\title[Localizations of inductively factored arrangements]
{Localizations of inductively factored arrangements}


\author[T. M\"oller]{Tilman M\"oller}
\address
{Fakult\"at f\"ur Mathematik,
Ruhr-Universit\"at Bochum,
D-44780 Bochum, Germany}
\email{tilman.moeller@rub.de}

\author[G. R\"ohrle]{Gerhard R\"ohrle}
\email{gerhard.roehrle@rub.de}


\allowdisplaybreaks

\begin{abstract}
We show that the class of
inductively factored arrangements
is closed under taking localizations.
We illustrate the usefulness of
this with an application.
\end{abstract}

\maketitle



\section{Introduction}

The notion of a nice arrangement is due to Terao \cite{terao:factored}.
This class generalizes the class of supersolvable arrangements, 
\cite{orliksolomonterao:hyperplanes}
(cf.\ \cite[Thm.\ 3.81]{orlikterao:arrangements}).
There is an inductive version of this class, 
so called inductively factored arrangements,
due to Jambu and Paris \cite{jambuparis:factored}, 
see Definition \ref{def:indfactored}.
This inductive class (properly) contains the class of supersolvable arrangements and 
is (properly) contained in the class of inductively free arrangements,
see \cite[Rem.\ 3.33]{hogeroehrle:factored}.

For an overview on properties of nice and inductively factored arrangements, 
and for their connection with the underlying Orlik-Solomon algebra, 
see \cite[\S 3]{orlikterao:arrangements}, \cite{jambuparis:factored}, 
and \cite{hogeroehrle:factored}.
In \cite{hogeroehrle:factored}, 
Hoge and the second author proved an 
addition-deletion theorem for nice arrangements, 
see Theorem \ref{thm:add-del-factored} below. 
This is an analogue of Terao's celebrated
Addition-Deletion Theorem \ref{thm:add-del} for 
free arrangements for the class of 
nice arrangements.

The class of free arrangements is known to be closed under 
taking localizations, \cite[Thm.\ 4.37]{orlikterao:arrangements}.
It is also known that this property restricts to 
various stronger notions of freeness, 
see \cite[Thm.\ 1.1]{hogeroehrleschauenburg:localizations}.
It is clear that the class of nice arrangements 
also satisfies this property, see Remark \ref{rem:factored}(ii) below.
Therefore, it is natural to investigate this question 
for the stronger property of inductively factored arrangements 
as well.
Here is the main result of our note.

\begin{theorem}
\label{thm:main}
The class of inductively factored arrangements is closed under
taking localizations.
\end{theorem}

Theorem \ref{thm:main} readily extends to the class of 
hereditarily inductively factored arrangements, see
Remark \ref{rem:heredindfac}.
Also, we give a short example to show 
the utility of such a result.

\section{Recollections and Preliminaries}
\label{sect:prelims}

\subsection{Hyperplane arrangements}
\label{ssect:arrangements}
Let $\BBK$ be a field and let
$V = \BBK^\ell$ be an $\ell$-dimensional $\BBK$-vector space.
A \emph{hyperplane arrangement} $\CA$ in $V$ 
is a finite collection of hyperplanes in $V$.
We also use the term $\ell$-arrangement for $\CA$. 

The \emph{lattice} $L(\CA)$ of $\CA$ is the set of subspaces of $V$ of
the form $H_1\cap \dotsm \cap H_i$ where $\{ H_1, \ldots, H_i\}$ is a subset
of $\CA$. 
For $X \in L(\CA)$, we have two associated arrangements, 
firstly
$\CA_X :=\{H \in \CA \mid X \subseteq H\} \subseteq \CA$,
the \emph{localization of $\CA$ at $X$}, 
and secondly, 
the \emph{restriction of $\CA$ to $X$}, $(\CA^X,X)$, where 
$\CA^X := \{ X \cap H \mid H \in \CA \setminus \CA_X\}$.
Note that $V$ belongs to $L(\CA)$
as the intersection of the empty 
collection of hyperplanes and $\CA^V = \CA$. 
The lattice $L(\CA)$ is a partially ordered set by reverse inclusion:
$X \le Y$ provided $Y \subseteq X$ for $X,Y \in L(\CA)$.

If $0 \in H$ for each $H$ in $\CA$, then 
$\CA$ is called \emph{central}.
If $\CA$ is central, then the \emph{center} 
$T_\CA := \cap_{H \in \CA} H$ of $\CA$ is the unique
maximal element in $L(\CA)$  with respect
to the partial order.
We have a \emph{rank} function on $L(\CA)$: $r(X) := \codim_V(X)$.
The \emph{rank} $r := r(\CA)$ of $\CA$ 
is the rank of a maximal element in $L(\CA)$.
Throughout, we only consider central arrangements.

More generally, for $U$ an arbitrary subspace of $V$, we can define  
$\CA_U :=\{H \in \CA \mid U \subseteq H\} \subseteq \CA$, the 
\emph{localization of $\CA$ at $U$}, 
and 
$\CA^U := \{ U \cap H \mid H \in \CA \setminus \CA_U\}$,
a subarrangement in $U$.

\subsection{Free hyperplane arrangements}
\label{ssect:free}
Let $S = S(V^*)$ be the symmetric algebra of the dual space $V^*$ of $V$.
Let $\Der(S)$ be the $S$-module of $\BBK$-derivations of $S$.
Since $S$ is graded, 
$\Der(S)$ is a graded $S$-module.

Let $\CA$ be an arrangement in $V$. 
Then for $H \in \CA$ we fix $\alpha_H \in V^*$ with
$H = \ker \alpha_H$.
The \emph{defining polynomial} $Q(\CA)$ of $\CA$ is given by 
$Q(\CA) := \prod_{H \in \CA} \alpha_H \in S$.
The \emph{module of $\CA$-derivations} of $\CA$ is 
defined by 
\[
D(\CA) := \{\theta \in \Der(S) \mid \theta(Q(\CA)) \in Q(\CA) S\} .
\]
We say that $\CA$ is \emph{free} if 
$D(\CA)$ is a free $S$-module, cf.\ \cite[\S 4]{orlikterao:arrangements}.

If $\CA$ is a free arrangement, then the $S$-module
$D(\CA)$ admits a basis of $n$ homogeneous derivations, 
say $\theta_1, \ldots, \theta_n$, \cite[Prop.\ 4.18]{orlikterao:arrangements}.
While the $\theta_i$'s are not unique, their polynomial 
degrees $\pdeg \theta_i$ 
are unique (up to ordering). This multiset is the set of 
\emph{exponents} of the free arrangement $\CA$
and is denoted by $\exp \CA$.

Terao's celebrated \emph{Addition-Deletion Theorem} 
which we recall next
plays a 
pivotal role in the study of free arrangements, 
\cite[\S 4]{orlikterao:arrangements}.
For $\CA$ non-empty, 
let $H_0 \in \CA$.
Define $\CA' := \CA \setminus\{ H_0\}$,
and $\CA'' := \CA^{H_0} = \{ H_0 \cap H \mid H \in \CA'\}$,
the restriction of $\CA$ to $H_0$.
Then $(\CA, \CA', \CA'')$ is a \emph{triple} of arrangements,
\cite[Def.\ 1.14]{orlikterao:arrangements}. 

\begin{theorem}[\cite{terao:freeI}]
\label{thm:add-del}
Suppose that $\CA$ is a non-empty $\ell$-arrangement.
Let  $(\CA, \CA', \CA'')$ be a triple of arrangements. Then any 
two of the following statements imply the third:
\begin{itemize}
\item[(i)] $\CA$ is free with $\exp \CA = \{ b_1, \ldots , b_{\ell -1}, b_\ell\}$;
\item[(ii)] $\CA'$ is free with $\exp \CA' = \{ b_1, \ldots , b_{\ell -1}, b_\ell-1\}$;
\item[(iii)] $\CA''$ is free with $\exp \CA'' = \{ b_1, \ldots , b_{\ell -1}\}$.
\end{itemize}
\end{theorem}
There are various stronger notions of freeness
which we discuss in the following subsections.

\subsection{Inductively free arrangements}
\label{ssect:indfree}

Theorem \ref{thm:add-del}
motivates the notion of 
\emph{inductively free} arrangements,  see 
\cite{terao:freeI} or 
\cite[Def.\ 4.53]{orlikterao:arrangements}.

\begin{defn}
\label{def:indfree}
The class $\CIF$ of \emph{inductively free} arrangements 
is the smallest class of arrangements subject to
\begin{itemize}
\item[(i)] $\Phi_\ell \in \CIF$ for each $\ell \ge 0$;
\item[(ii)] if there exists a hyperplane $H_0 \in \CA$ such that both
$\CA'$ and $\CA''$ belong to $\CIF$, and $\exp \CA '' \subseteq \exp \CA'$, 
then $\CA$ also belongs to $\CIF$.
\end{itemize}
\end{defn}

Free arrangements are closed with respect to taking 
localizations, 
cf.~\cite[Thm.\ 4.37]{orlikterao:arrangements}. 
This also holds for the class $\CIF$.

\begin{theorem}
[{\cite[Thm.\ 1.1]{hogeroehrleschauenburg:localizations}}]
\label{thm:indfreelocal}
If $\CA$ is inductively free,
then so is $\CA_U$ for every subspace $U$ in $V$.
\end{theorem}

\subsection{Nice and inductively factored arrangements}
\label{ssect:factored}

The notion of a \emph{nice} or \emph{factored} 
arrangement goes back to Terao \cite{terao:factored}.
It generalizes the concept of a supersolvable arrangement, see
\cite[Thm.\ 5.3]{orliksolomonterao:hyperplanes} and 
\cite[Prop.\ 2.67, Thm.\ 3.81]{orlikterao:arrangements}.
Terao's main motivation was to give a 
general combinatorial framework to 
deduce factorizations of the underlying Orlik-Solomon algebra,
see also \cite[\S 3.3]{orlikterao:arrangements}.
We recall the relevant notions  
from \cite{terao:factored}
(cf.\  \cite[\S 2.3]{orlikterao:arrangements}):

\begin{defn}
\label{def:factored}
Let $\pi = (\pi_1, \ldots , \pi_s)$ be a partition of $\CA$.
\begin{itemize}
\item[(a)]
$\pi$ is called \emph{independent}, provided 
for any choice $H_i \in \pi_i$ for $1 \le i \le s$,
the resulting $s$ hyperplanes are linearly independent, i.e.\
$r(H_1 \cap \ldots \cap H_s) = s$.
\item[(b)]
Let $X \in L(\CA)$.
The \emph{induced partition} $\pi_X$ of $\CA_X$ is given by the non-empty 
blocks of the form $\pi_i \cap \CA_X$.
\item[(c)]
$\pi$ is
\emph{nice} for $\CA$ or a \emph{factorization} of $\CA$  provided 
\begin{itemize}
\item[(i)] $\pi$ is independent, and 
\item[(ii)] for each $X \in L(\CA) \setminus \{V\}$, the induced partition $\pi_X$ admits a block 
which is a singleton. 
\end{itemize}
\end{itemize}
If $\CA$ admits a factorization, then we also say that $\CA$ is \emph{factored} or \emph{nice}.
\end{defn}

\begin{remark}
\label{rem:factored}
The class of nice arrangements is closed under taking localizations.
For, if $\CA$ is non-empty and   
$\pi$ is a nice partition of $\CA$, then the non-empty parts of the 
induced partition $\pi_X$ form a nice partition of $\CA_X$
for each $X \in L(\CA)\setminus\{V\}$;
cf.~the proof of \cite[Cor.\ 2.11]{terao:factored}.
\end{remark}

Following Jambu and Paris 
\cite{jambuparis:factored}, 
we introduce further notation.
Suppose $\CA$ is not empty. 
Let $\pi = (\pi_1, \ldots, \pi_s)$ be a partition of $\CA$.
Let $H_0 \in \pi_1$ and 
let $(\CA, \CA', \CA'')$ be the triple associated with $H_0$. 
Then $\pi$ induces a 
partition 
$\pi'$ of 
$\CA'$, i.e.\ the non-empty 
subsets $\pi_i \cap \CA'$.
Note that since $H_0 \in \pi_1$, we have
$\pi_i \cap \CA' = \pi_i$ 
for $i = 2, \ldots, s$. 
Also, associated with $\pi$ and $H_0$, we define 
the \emph{restriction map}
\[
\R := \R_{\pi,H_0} : \CA \setminus \pi_1 \to \CA''\ \text{ given by } \ H \mapsto H \cap H_0
\]
and set 
\[
\pi_i'' := \R(\pi_i) = \{H \cap H_0 \mid H \in \pi_i\} \
\text{ for }\  2 \le i \le s.
\]
In general, $\R$ need not be surjective nor injective.
However, since we are only concerned with cases when 
$\pi'' = (\pi_2'', \ldots, \pi_s'')$ is a
partition of $\CA''$,  
$\R$ has to be onto and 
$\R(\pi_i) \cap \R(\pi_j) = \varnothing$ for $i \ne j$.

The following 
analogue of Terao's 
Addition-Deletion Theorem \ref{thm:add-del} for 
free arrangements for the class of 
nice arrangements is proved in 
\cite[Thm.\ 3.5]{hogeroehrle:factored}.

\begin{theorem}
\label{thm:add-del-factored}
Suppose that $\CA \ne \Phi_\ell$.
Let $\pi = (\pi_1, \ldots, \pi_s)$ be a  partition  of $\CA$.
Let $H_0 \in \pi_1$ and 
let $(\CA, \CA', \CA'')$ be the triple associated with $H_0$. 
Then any two of the following statements imply the third:
\begin{itemize}
\item[(i)] $\pi$ is nice for $\CA$;
\item[(ii)] $\pi'$ is nice for $\CA'$;
\item[(iii)] $\R: \CA \setminus \pi_1 \to \CA''$ 
is bijective and $\pi''$ is nice for $\CA''$.
\end{itemize}
\end{theorem}

Note the bijectivity condition on $\R$ 
in Theorem \ref{thm:add-del-factored}
is necessary, 
cf.~\cite[Ex.\ 3.3]{hogeroehrle:factored}.
Theorem \ref{thm:add-del-factored} 
motivates
the following stronger notion of factorization, 
cf.\ \cite{jambuparis:factored}, \cite[Def.\ 3.8]{hogeroehrle:factored}.

\begin{defn} 
\label{def:indfactored}
The class $\CIFAC$ of \emph{inductively factored} arrangements 
is the smallest class of pairs $(\CA, \pi)$ of 
arrangements $\CA$ together with a partition $\pi$
subject to
\begin{itemize}
\item[(i)] 
$(\Phi_\ell, (\varnothing)) \in \CIFAC$ for each $\ell \ge 0$;
\item[(ii)] 
if there exists a partition $\pi$ of $\CA$ 
and a hyperplane $H_0 \in \pi_1$ such that 
for the triple $(\CA, \CA', \CA'')$ associated with $H_0$ 
the restriction map $\R = \R_{\pi, H_0} : \CA \setminus \pi_1 \to \CA''$ 
is bijective and for the induced partitions $\pi'$ of $\CA'$ and 
$\pi''$ of $\CA''$ 
both $(\CA', \pi')$ and $(\CA'', \pi'')$ belong to $\CIFAC$, 
then $(\CA, \pi)$ also belongs to $\CIFAC$.
\end{itemize}
If $(\CA, \pi)$ is in $\CIFAC$, then we say that
$\CA$ is \emph{inductively factored with respect to $\pi$}, or else
that $\pi$ is an \emph{inductive factorization} of $\CA$. 
Sometimes we simply say $\CA$ is \emph{inductively factored} without 
reference to a specific inductive factorization of $\CA$.
\end{defn}

\begin{remark}
\label{rem:indtable}
If $\pi$ is an inductive factorization of $\CA$, then  
there exists an \emph{induction of factorizations} 
by means of Theorem \ref{thm:add-del-factored} as follows.
This procedure amounts to 
choosing a total order on $\CA$, say 
$\CA = \{H_1, \ldots, H_n\}$, 
so that each of the pairs 
$(\CA_0 = \Phi_\ell, \varnothing)$,
$(\CA_i := \{H_1, \ldots, H_i\}, \pi_i:= \pi|_{\CA_i})$, and 
$(\CA_i'' :=\CA_i^{H_i}, \pi_i'')$
for each $1 \le i \le n$, 
belongs to $\CIFAC$ 
see \cite[Rem.\ 3.16]{hogeroehrle:factored}.
\end{remark}

The connection with the previous notions is as follows.

\begin{proposition}
[{\cite[Prop.\ 3.11]{hogeroehrle:factored}}]
\label{prop:superindfactored}
If $\CA$ is supersolvable, then $\CA$ is inductively factored.
\end{proposition}

\begin{proposition}
[{\cite[Prop.\ 2.2]{jambuparis:factored}, 
\cite[Prop.\ 3.14]{hogeroehrle:factored}}]
\label{prop:indfactoredindfree}
Let $\pi = (\pi_1, \ldots, \pi_r)$ be an inductive factorization of $\CA$. 
Then $\CA$ is inductively free with 
$\exp \CA = \{0^{\ell-r}, |\pi_1|, \ldots, |\pi_r|\}$.
\end{proposition}

\begin{defn}
\label{def:heredindfactored}
We say that $\CA$ 
is \emph{hereditarily inductively factored}
provided $\CA^Y$ is inductively factored for every $Y \in L(\CA)$.
\end{defn}

\section{Proof of Theorem \ref{thm:main}}
\label{sec:proof}

We readily reduce to the case where 
we localize with respect to a space
$X$ belonging to the 
intersection lattice of $\CA$. 
For, letting $X = \cap_{H \in \CA_U} H \in L(\CA)$,
we have $\CA_X = \CA_U$.

We are going to show that if $\pi$ is 
an inductive factorization of $\CA$, then 
the restriction $\pi_X$ of $\pi$ to $\CA_X$ is 
an inductive factorization of the latter.
We argue by induction on the rank $r(\CA)$.
If $r(\CA) \le 3$, then 
$r(\CA_X) \le 2$ for $X \ne T_\CA$, so
the result follows from the proof of 
Proposition \ref{prop:superindfactored}
(and the fact that $V < H < X = T(\CA_X)$ is a
maximal chain of modular elements in $L(\CA_X)$
for every $H \in L(\CA_X)$).

So suppose $\CA$ is inductively factored of rank $r > 3$
and that the statement above holds for 
all inductively factored arrangements of rank less than $r$.
Let $\pi$ be an inductive factorization of $\CA$.
Let $\{H_1, \ldots H_n\}$ be the total order on 
$\CA$ such that for $\CA_i := \{H_1, \ldots H_i\}$
the induced partition $\pi_i := \pi|_{\CA_i}$ is an inductive factorization of 
$\CA_i$ for $i = 1, \ldots, n$, see Remark \ref{rem:indtable}.
Consider the sequence of 
inductive factorizations
\begin{equation}
\label{eq:seq0}
(\CA_1,\pi_1),  (\CA_2,\pi_2),  \ldots, 
(\CA_n,\pi_n) = (\CA, \pi).
\end{equation}

Then, for $i = 1, \ldots, n$, we have
\begin{equation}
\label{eq:seq1}
\CA_X \cap \CA_i = (\CA_i)_X.
\end{equation}
For $H \in \CA_X \cap \CA_i$, we have $H \le X$, and so by \eqref{eq:seq1},
for $i = 1, \ldots, n$,  
\begin{equation}
\label{eq:seq2}
\left(\CA_X \cap \CA_i\right)^{H} = \left((\CA_i)_X\right)^{H} 
= \left(\CA_i^{H}\right)_X.
\end{equation}

Consequently, localizing each member $(\CA_i,\pi_i)$ of 
the sequence \eqref{eq:seq0} at $X$, removing redundant 
terms if necessary and reindexing the resulting distinct arrangements, 
we obtain the following sequence of 
subarrangements of 
$\CA_X$, 
\begin{equation}
\label{eq:seq3}
\CA_{1,X} \subsetneq \CA_{2,X}
\subsetneq  \ldots \subsetneq 
\CA_{m,X} = \CA_X,
\end{equation}
where $\CA_{i,X}$ is short for $(\CA_i)_X$. 
In particular, $|\CA_{i,X}| = i$ and $m \le n$.
Thus we obtain the following sequence of subarrangements of 
$\CA_X$ along with induced partitions:
\begin{equation}
\label{eq:seq4}
(\CA_{1,X}, \pi_{1,X}), (\CA_{2,X}, \pi_{2,X}), \ldots,
(\CA_{m,X}, \pi_{m,X}) = (\CA_X, \pi_X),
\end{equation}
where $\pi_{i,X}$ is the induced partition of $\pi_i$ on $\CA_{i,X}$,
i.e.\ $\pi_{i,X} := \pi|_{\CA_{i,X}}$.
We claim that \eqref{eq:seq4} 
is an inductive chain of factorizations of
$\CA_X$, so that
$\pi_X$ is an inductive factorization of $\CA_X$.

Now let $H_i \in \CA_X \cap \CA_i = \CA_{i,X}$ be the relevant 
hyperplane in the $i$th step in the sequence \eqref{eq:seq3}.
Let $(\CA_{i,X}, \CA_{i,X}', \CA_{i,X}'')$
be the triple with respect to $H_i$.
Thus, by the constructions of the chains 
in \eqref{eq:seq3} and  \eqref{eq:seq4}, 
we have for each $2 \le i \le m$
\begin{equation}
\label{eq:seq8}
\CA_{i,X}' =  \CA_{i-1, X}
\quad \text{and} \quad \pi_{i,X}' = \pi_{i-1,X}. 
\end{equation}
Without loss, suppose $H_i$ belongs to 
the first part $(\pi_{i,X})_1$ of $\pi_{i,X}$ for each $i$.
Thanks to \eqref{eq:seq8}, 
Remark \ref{rem:factored} and 
Theorem  \ref{thm:add-del-factored}
the corresponding restriction map 
$\R : \CA_{i,X} \setminus (\pi_{i,X})_1 \to \CA_{i,X}''$
is bijective for each $i$.

Since $(\CA_i'', \pi_i'')$ belongs to $\CIFAC$, it follows 
by induction on the rank that each localization 
\[
(\CA_{i,X}'', \pi_{i,X}'') :=
((\CA_i'')_X, (\pi_i'')_{X})
\] 
also belongs to $\CIFAC$, for each $i < m$,
where we used  \eqref{eq:seq4} and set 
$\pi_{i,X}'' := (\pi_i'')_{X}$.
 


Finally, since $\CA_{1,X}$ is of rank $1$, 
$(\CA_{1,X}, \pi_{1,X})$ belongs to $\CIFAC$.
Thus, because $(\CA_{i,X}'', \pi_{i,X}'')$ 
belongs to $\CIFAC$ 
and the corresponding restriction map
$\R  : \CA_{i,X} \setminus (\pi_{i,X})_1 \to \CA_{i,X}''$
is bijective for each $i < m$, 
it follows from \eqref{eq:seq8} and 
a repeated application of the
addition part of Theorem \ref{thm:add-del-factored}
that also 
$(\CA_{m,X}, \pi_{m,X}) = (\CA_{X}, \pi_{X})$
 belongs to $\CIFAC$, as desired.

\begin{remark}
\label{rem:heredindfac}
Theorem \ref{thm:main} readily extends to hereditarily inductively factored 
arrangements. For, let $\CA$ be hereditarily inductively factored and let 
$Y \le X$ in $L(\CA)$.
Then, since $\CA^Y$ is inductively factored, so is $(\CA^Y)_X$, 
by Theorem \ref{thm:main}.
Finally, since $(\CA_X)^Y = (\CA^Y)_X$,
it follows that $(\CA_X)^Y $ is inductively factored.
\end{remark}

The following example shows the utility of the results above.

\begin{example}
\label{ex:akl}
Let $V = \BBC^\ell$ be an $\ell$-dimensional $\BBC$-vector space.
Orlik and Solomon defined intermediate 
arrangements $\CA^k_\ell(r)$ in 
\cite[\S 2]{orliksolomon:unitaryreflectiongroups}
(cf.\ \cite[\S 6.4]{orlikterao:arrangements}) which
interpolate between the
reflection arrangements 
$\CA(G(r,1,\ell))$ and $\CA(G(r,r,\ell))$ of 
the complex reflection groups
$G(r,1,\ell)$ and $G(r,r,\ell)$. 
For  $\ell, r \geq 2$ and $0 \leq k \leq \ell$, the defining polynomial of
$\CA^k_\ell(r)$ is 
$$
Q(\CA^k_\ell(r)) = x_1 \cdots x_k\prod\limits_{\substack{1 \leq i < j \leq \ell\\ 0 \leq n < r}}(x_i - \zeta^nx_j),
$$
where $\zeta$ is a primitive $r$th root of unity,
so that 
$\CA^\ell_\ell(r) = \CA(G(r,1,\ell))$ and 
$\CA^0_\ell(r) = \CA(G(r,r,\ell))$.
Note that for $1 < k < \ell$, 
$\CA^k_\ell(r)$ is not a reflection arrangement.

Each of these arrangements is known to be free,
cf.~\cite[Prop.~6.85]{orlikterao:arrangements}.
The supersolvable and inductively free cases among them 
were classified in \cite{amendhogeroehrle:super},
and \cite{amendhogeroehrle:indfree}, respectively. 

If $k \in \{\ell-1, \ell\}$, then $\CA^k_\ell(r)$ is supersolvable, 
by \cite[Thm.\ 1.3]{amendhogeroehrle:super},
and so $\CA^k_\ell(r)$ is inductively factored,
by Proposition \ref{prop:superindfactored}.
Let $\ell \ge 4$. We claim that 
$\CA^k_\ell(r)$ is not nice for $0 \le k \le \ell-4$ and moreover
$\CA^k_\ell(r)$ is not inductively factored for $0 \le k \le \ell-3$.

For $k = 0$, this follows from \cite[Thm.\ 1.3]{hogeroehrle:nice}.
So let $1 \le k \le \ell-3$ and set $\CA = \CA^k_\ell(r)$. Define
\[
X := \bigcap\limits_{\substack{k+1 \leq i < j \leq \ell\\ 0 \leq n < r}} \ker(x_i - \zeta^nx_j).
\]
Then one checks that 
\[
\CA_X \cong \CA^0_{\ell-k}(r) = \CA(G(r,r,\ell-k)).
\]
For $1 \le k \le \ell-4$, it follows from \cite[Thm.\ 1.3]{hogeroehrle:nice} 
that $\CA(G(r,r,\ell-k))$ is not nice. Consequently, neither is 
$\CA^k_\ell(r)$, by Remark \ref{rem:factored}.
For $k = \ell-3$, we have 
$\CA_X \cong \CA(G(r,r,3))$.
By \cite[Cor.\ 1.4]{hogeroehrle:nice},
the latter is not inductively factored, 
thus  neither is $\CA^{\ell-3}_\ell(r)$,
thanks to Theorem \ref{thm:main}.
\end{example}

\medskip
{\bf Acknowledgments}:
We acknowledge 
support from the DFG-priority program 
SPP1489 ``Algorithmic and Experimental Methods in
Algebra, Geometry, and Number Theory''.


\bigskip

\bibliographystyle{amsalpha}

\newcommand{\etalchar}[1]{$^{#1}$}
\providecommand{\bysame}{\leavevmode\hbox to3em{\hrulefill}\thinspace}
\providecommand{\MR}{\relax\ifhmode\unskip\space\fi MR }
\providecommand{\MRhref}[2]{%
  \href{http://www.ams.org/mathscinet-getitem?mr=#1}{#2} }
\providecommand{\href}[2]{#2}


\end{document}